\documentclass[11pt]{article}

\usepackage{amsthm}

\long\def\eatit#1{}

\newtheorem{thm}{Theorem}[subsection]

\newtheorem{lem}[thm]{Lemma}
\newtheorem{cor}[thm]{Corollary}

\newtheorem{Rmk}[thm]{Remark}

\def\pr#1{\hbox{{\bf P}${}^{#1}$}}
\def\sc#1{\hbox{$\varepsilon(#1)$}}

\begin{document}

\title{Discrete Behavior of Seshadri Constants on Surfaces}

\author{Brian harbourne\\
Department of Mathematics\\
University of Nebraska-Lincoln\\
Lincoln, NE 68588-0130\\
USA\\
WEB: http://www.math.unl.edu/$\scriptstyle\sim$bharbour/ \\
email: bharbour@math.unl.edu\\
\and
Joaquim Ro\'e\\
Departament de Matem\`atiques\\   
Universitat Aut\`onoma de Barcelona\\
08193 Bellaterra (Barcelona)\\
Spain\\
email: jroe@mat.uab.es}

\date{June 7, 2007}

\maketitle

\begin{abstract} Working over ${\bf C}$,
we show that, apart possibly from a unique limit point, the 
possible values of multi-point Seshadri constants 
for general points on
smooth projective surfaces form a discrete set. In addition to
its theoretical interest, this result is of practical value,
which we demonstrate by giving
significantly improved explicit lower bounds
for Seshadri constants on \pr2 and new 
results about ample divisors on blow ups of \pr2 at general points.
\end{abstract}

\thanks{Acknowledgments: We would like to thank T.\ Szemberg for making  
his work \cite{refS} available to us.
Harbourne would also like to thank the NSA and the NSF 
for their support, and Ro\'e
would like to thank the support of the {\it Programa Ram\'on y
Cajal} of the spanish MCyT, and of the projects CAICYT BFM2002-01240,
2000SGR-00028 and EAGER. \newline
\hbox to.2in{\hfil} {\it 2000 Mathematics Subject Classification}. Primary 14C20; Secondary 14J99.\newline
\hbox to.2in{\hfil} {\it Key words and phrases}. Multi-point Seshadri constants
on surfaces, ample divisors.}

\section{Introduction}\label{intro}

The situation often arises that one has a birational
morphism of smooth projective varieties $\pi: Y\to X$, where $X$ is well understood 
and one wants to understand $Y$. For example, even if one
knows precisely which divisors on $X$ are ample, or nef,
it is often a difficult problem to determine the same for $Y$.
The problem of determining ampleness or nefness on $Y$ is closely related to
the problem of computing multi-point Seshadri constants on $X$.

Even in the case that $X$ is a surface, it is quite hard to compute Seshadri 
constants exactly. Our approach instead is to study what values are possible.
Of course, the more one knows about a surface $X$ the more one would hope to
be able to restrict what is possible. What has not been previously recognized is
that easily obtained information about $X$ already puts a lot of structure on 
the set of possible values of Seshadri constants: 
if the blown up points are general, the set of possible
values is, apart possibly from a unique limit point, a discrete set.
This has significant consequences for determining Seshadri constants on surfaces;
one consequence, for example, is our Theorem \ref{thetheorem}, 
which establishes a framework for computing arbitrarily accurate lower 
bounds for multi-point Seshadri constants.
Although we do not focus on implementing this framework here
(for a detailed consideration of algorithmic concerns,
see the unpublished posting \cite{refoldpreprint}),
we do demonstrate what our methods can achieve with results easily at hand
by giving significant improvements to previously known lower bounds
for multi-point homogeneous Seshadri constants on \pr2 (Corollary \ref{corRem}), 
and we determine ampleness for many new cases
on blow ups $Y$ of \pr2 at general points (Corollary \ref{simplifiedampCor}). 

\subsection{Seshadri constants}\label{screview}

Let $X$ be a smooth projective variety of dimension $N>1$, and let $L$ be a nef
divisor class (i.e., $L^r\cdot Z\ge0$ for 
every effective $r$-cycle $Z$ on $X$).
Given a positive integer $n$ and a nonzero real vector
${\bf \ell}=(l_1,\cdots,l_n)$ with each $l_i\ge0$, the multi-point
Seshadri constant for ${\bf \ell}$ and points $p_1, \dots, p_n$
of $X$ is the real number
$$\varepsilon(X,L,l_1p_1, \dots, l_np_n)= {\hbox{inf}
\left\{{L\cdot C}\over{\Sigma_{i=1}^n 
l_i\hbox{mult}_{p_i}C}\right\}},$$ 
where the infimum is taken with respect to 
all curves $C$ through at least 
one of the points. For the one-point and the multi-point
homogeneous case (in which 
$l_i=1$ for all $i$ and which most previous work has focused on), 
see \cite{refDe} or \cite{refSS}. 
We also take $\varepsilon(X,L,n,{\bf \ell})$ to be defined as 
${\hbox{sup}\{\varepsilon(X,L,l_1p_1, \dots, l_np_n)\}},$ where
the supremum is taken with respect to 
all choices of $n$ distinct points $p_i$ of $X$.
For the homogeneous case, we write simply $\varepsilon(X,L,n)$
in place of $\varepsilon(X,L,n,(1,\ldots,1))$.
Since the homogeneous case where $X=\pr2$ and $L$ is
the class of a line is of particular interest, we will denote \sc{\pr2,L,n} simply by \sc{n}.

It is well known and not difficult to prove that 
$\sc{X,L,p_1, \dots, p_n}\le \root{N}\of{L^{N}/n}$,
but lower bounds are much more challenging (see 
\cite{refN}, \cite{refKu} and \cite{refXb}).
It is not
hard to see that  $\varepsilon(X,L,n)=\varepsilon(X,L,p_1, \dots, p_n)$
for {\it very general} points $p_1, \dots, p_n$ 
(i.e., in the intersection
of countably many Zariski-open and dense subsets of $X^n$), 
   although some results (see \cite{refO}, \cite{refSS}) suggest that the
   equality might hold in fact for {\it general} points (i.e., in a
   Zariski-open subset of $X^n$). 
When $L$ is a big (i.e., $L^2>0$) and nef divisor on a surface $X$, 
our Theorem \ref{thetheorem} gives lower bounds for $\varepsilon(X,L,n)$
which in fact hold for $\varepsilon(X,L,p_1, \dots, p_n)$ 
for general points $p_i$. 

Two methods have been used to give lower bounds on \sc{X,L,n}
for surfaces $X$. One involves explicit
constructions of nef divisors, the other involves ruling out the 
existence of certain putative reduced irreducible curves
of negative self-intersection (so-called $L$-abnormal curves). Both methods,
which work also in the non-homogeneous case, depend
on looking at the surface $Y$ obtained from $X$ by the morphism
$\pi:Y\rightarrow X$ blowing up distinct points $p_i\in X$, $1\le i\le n$.
If $E_i$ is the divisor
class of the exceptional curve $\pi^{-1}(p_i)$, then clearly
$\varepsilon(X,L,l_1p_1, \dots, l_np_n)$ is the largest $t$ such that
$F_t=\pi^*L-t(l_1E_1+\cdots+l_nE_n)$ is nef, hence
$\sc{X,L,n,{\bf \ell}}\ge t$ whenever $F_t=\pi^*L-t(l_1E_1+\cdots+l_nE_n)$ 
is a nef ${\bf R}$-divisor class 
    (i.e., a nef element of the 
    divisor class group with real coefficients). 

Alternatively (see Lemma \ref{introlem}), suppose each $l_i$ is rational and $t$,
$0\le t<\sqrt{L^2/{\bf \ell}^2}$, is rational, where
${\bf \ell}^2$ signifies the usual dot product.
Then $t\le\sc{X,L,n,{\bf \ell}}$ if and only if, for general points $p_i$
there are no reduced and irreducible curves $C\subset X$ such that
$F_t\cdot H<0$ where 
$H=\pi^*C-h_1E_1-\cdots-h_nE_n$ is the class of
the proper transform of $C$ (so $h_i$ is the multiplicity of $C$ at
$p_i$); note that $F_t\cdot H<0$ is equivalent to
$(L\cdot C)/(l_1h_1+\cdots+l_nh_n) < t$).
In the homogeneous case we call such a curve $C$ an $L$-{\it abnormal\/} curve 
(or simply {\it abnormal\/} if $L$ is understood),
following Nagata \cite{refN}, 
    who, in case ${\bf \ell}=(1,\ldots,1)$ and
    $L$ is a line in $X=\pr2$, called any such 
    curve $C$ an abnormal curve (also referred to as 
    submaximal in \cite{refBa} and \cite{refSS}). 
Moreover, if $\hbox{Pic}(X)/\hbox{$\sim$}$, where $\sim$ denotes numerical
equivalence, is cyclic (as is the case for $X=\pr2$), 
then for any such $C$ we have $\sc{X,L,n} = (L\cdot
C)/(h_1+\cdots+h_n)$ by Lemma \ref{introlemb}.
    (For \pr2, Nagata also found all curves 
    abnormal for each $n<10$, showed
    no curve is abnormal for $n$ when $n$ is a square and conjectured
    there are no abnormal curves for $n\ge10$.)

    So, to exemplify the first method, if for some choice of distinct
    points $p_i$ one finds positive integers $d$ and $t$ such that
    $d\pi^*L-t(l_1E_1+\cdots+l_nE_n)$ is  nef, it follows that
    $\sc{X,L,n,{\bf \ell}}\ge d/t$.  This basic idea is used in
    \cite{refBi} (for $X=\pr2$) and \cite{refH} (for surfaces
    generally) to obtain bounds of the form $\sc{X,L,n}\ge
    (\sqrt{L^2/n})\sqrt{1-1/f(n)}$ where $f(n)$, for some values of
    $n$, is a quadratic function of $n$.  
        Note that the bound $\sc{n}\ge(1/\sqrt{n})(\sqrt{1-1/f(n)})$
        is equivalent to the inequality ${\cal R}_n(L) \le 1/f(n)$ of
        \cite{refBi}, where ${\cal R}_n(L)$ is  what is called in
        \cite{refBi} the {\it $n$-th remainder\/}  of the divisor
        class $L$.
    Alternatively,   
to exemplify the second method,
suppose one is given $F_t=\pi^*L-t(l_1E_1+\cdots+l_nE_n)$.
One then constructs a set $o_n(F_t)$ of values
which one somehow can show contains  
$(\pi^*L\cdot D)/(-(l_1E_1+\cdots+l_nE_n)\cdot D)$ 
for every effective, reduced, irreducible divisor $D$
on $Y$ with $F_t\cdot D<0$, if any. (We show how to obtain a 
specific such set $o_n(F_t)$ after Lemma \ref{testLem}.) 
For as many values $v\in o_n(F_t)$ as possible,
one attempts to show that there is no such $D$ for which
$v=(\pi^*L\cdot D)/(-(l_1E_1+\cdots+l_nE_n)\cdot D)$. If $c$ is the infimum 
of the remaining values in $o_n(F_t)$, then we conclude that
$F_c$ is nef and hence that $c\le \sc{X,L,n,{\bf \ell}}$.
Thus the more values $v\in o_n(F_t)$ one can rule out, the better this bound
becomes. For the homogeneous case,  this is the basic idea used
implicitly in \cite{refX}, \cite{refSS}, \cite{refST} and \cite{refT},
with the latter obtaining the bound
$\sc{n}\ge(1/\sqrt{n})\sqrt{1-1/(12\,n+1)}$.

Given ${\bf \ell}$ and a big and nef $L$, we can, 
for each $c<\sqrt{L^2/{\bf \ell}^2}$, give a finite set 
$o_n(F_c)$ (see Theorem \ref{maindiscthm}) depending only on $\ell$, $c$, $L^2$
and the semigroup of $L$-degrees $\{C\cdot L : C$ is an effective
divisor$\}$ of curves. 
This shows the set of possible values
of $\varepsilon(X,L,n,{\bf \ell})$ is either finite or an increasing 
discrete sequence and, in the latter case, 
$\sqrt{L^2/{\bf \ell}^2}$ is its unique limit point.
I.e., apart from $\sqrt{L^2/{\bf \ell}^2}$, the set of possible values of $\varepsilon(X,L,n,{\bf \ell})$ is discrete.
This has a number of conceptual consequences. For example, 
if we write this increasing sequence as $o(n,L)_1 < o(n,L)_2 < \cdots$,
and if we were to show that $o(n,L)_i<\varepsilon(X,L,n)$, then in fact
it automatically follows that $o(n,L)_{i+1}\le\varepsilon(X,L,n)$. 
Moreover, to show $\varepsilon(X,L,n,{\bf \ell})\ge c$ for any $c<\sqrt{L^2/{\bf \ell}^2}$, 
there are only finitely many values of $\varepsilon(X,L,n,{\bf \ell})$ less than 
$c$ one must rule out. Moreover,
carrying this calculation out will either show that
$\varepsilon(X,L,n,{\bf \ell})\ge c$, or it will compute
$\varepsilon(X,L,n,{\bf \ell})$ exactly (by finding which value in $o_n(F_c)$
is the correct one).

Our general results about the existence of $o_n(F_c)$
with the structure as claimed above
are stated in Theorem \ref{maindiscthm} and proved in Subsection \ref{disc}.
Using refinements of these results which we obtain in Subsection \ref{appl}, 
we then prove Theorem \ref{thetheorem} (which shows how
theoretical results ruling out the existence of abnormal curves can be
converted into bounds on \sc{X,L,n}) and
Corollary \ref{corRem} (which gives lower bounds for $\sc{n}$ 
that for most values of $n$ are significantly
better than what was known previously). As another
application, we also obtain in Corollary \ref{simplifiedampCor} 
improved results on ample divisors
on blow ups of \pr2.

\subsection{Applications}\label{mrintro}

Our results involve
a related apparently simpler problem, that of the existence of
curves with a given sequence of 
multiplicities ${\bf m}=(m_1, \dots,
m_n)$ at given points $p_1$, \dots, $p_n \in X$. 
Let us denote by
$\alpha(X,L,{\bf m},p_1, \dots, p_n)$ (respectively, 
$\alpha_0(X,L,{\bf m},p_1, \dots, p_n)$) the least 
degree $L \cdot C$ of a curve $C$ (respectively,
irreducible curve)  
passing with multiplicity at least $m_i$ (respectively, exactly $m_i$)
through each  point  $p_i$. If the points are in general position in
$X$, we write simply $\alpha(X,L,{\bf m})$ and $\alpha_0(X,L,{\bf
m})$. When focusing on the case that $L$ is a line in $X=\pr{2}$, we will denote
$\alpha(\pr2,L,{\bf m})$ and $\alpha_0(\pr2,L,{\bf m})$ 
simply by $\alpha({\bf m})$ and $\alpha_0({\bf m})$. 
Given an integer $m$, we will denote the vector $(m,\ldots,m)$ with $r$
entries of $m$ by $m^{[r]}$. As a consequence of our 
results in Section \ref{submax}, we will prove the following:

\begin{thm}\label{thetheorem}
Let $X$ be a smooth projective surface,
$L$ a big and nef divisor, $n\ge 2$ an integer and $\mu \ge 1$ a real number.
\begin{itemize}
\item[(a)] If $\alpha(X,L,m^{[n]}) \ge m \sqrt{L^2(n-1/\mu)}$ for
every integer $1 \le m < \mu$,  
then  $$\sc{X,L,n} >\sqrt{L^2\over n}\ \sqrt{1-{1\over{(n-2)\mu}}}.$$
\item[(b)] If $\alpha_0(X,L,m^{[n]}) \ge m \sqrt{L^2(n-1/\mu)}$
for every integer $1 \le m < \mu$, and if
$$\alpha_0((m^{[n-1]},m+k)) \ge {{mn+k}\over{n}} \sqrt{L^2(n-1/\mu)}$$ 
for every integer $1 \le m < \mu/(n-1)$
and every integer $k$ with 
$$k^2 < (n/(n-1))\hbox{min}\,(m,m+k),$$
then $$\sc{X,L,n} \ge \sqrt{L^2\over n}\ \sqrt{1-{1\over{n\mu}}}.$$
\end{itemize}
\end{thm}

In order to apply the theorem, one just needs to know some values of
$\alpha$. 
Drawing on asymptotic results of Alexander and Hirschowitz, 
for example, it is
possible to give bounds on $\varepsilon$ for surfaces on which the Picard
group is generated by a single ample divisor. In fact, the main result
of \cite{refAH} already implies 
ampleness for certain divisors (and so bounds on $\varepsilon(X,L,n)$ for
some n); a suitable interpretation of Theorem \ref{thetheorem} 
yields the following corollary,
linking the mentioned asymptotic results to lower bounds
for Seshadri constants in the form  $(\sqrt{L^2/n}\,)\sqrt{1-1/f(n)}$,
analogous to what is known for \pr2. 

\begin{cor}\label{corAH}
Let $X$ be a surface on which the Picard
group is generated by a single ample divisor $L$, and let  ${\bf m}:{\bf N
\rightarrow N}$ be a map such that for every $m< {\bf m}(n)$
\begin{equation}
\label{ahbound}
\dim \Bigg|{{\alpha(X,L,m^{[n]})}\over{ L^2}}L\Bigg|\ge
n{{m(m+1)}\over 2}
\end{equation}
holds. Then there is an $n_0$ such that for $n\ge n_0$, 
$$\sc{X,L,n}\ge \sqrt{{L^2}\over {n}}\sqrt{1-{1\over
(n-2){\bf m}(n)}}\,.$$
Moreover, there exists such an ${\bf m}(n)$ with 
$\displaystyle\lim_{n\to \infty}{\bf m}(n)=\infty$ and hence
$\displaystyle\lim_{n\to \infty}n{\cal R}_n(L)=0$.
\end{cor}

We can give much more specific bounds
for \pr2. For instance, for $X=\pr2$
it is known that $\alpha(m^{[n]}) \ge m\sqrt{n}$ for $n\ge10$ and
$m\le \lfloor \sqrt{n}\rfloor(\lfloor \sqrt{n}\rfloor-3)/2$
(see the proof of Corollary 1.2(a) of \cite{refHR}), so we may apply 
Theorem \ref{thetheorem}(b) with
$\mu=1+\lfloor\sqrt{n}\rfloor(\lfloor \sqrt{n}\rfloor-3)/2$
whenever $1<1+\lfloor\sqrt{n}\rfloor(\lfloor \sqrt{n}\rfloor-3)/2$,
so for $n\ge 16$.
(Note that the hypotheses involving $k\ne 0$ are vacuous when
$\mu/(n-1)<1$.)
On the other hand, results of \cite{refCCMO} imply
that $\alpha(m^{[n]}) \ge m\sqrt{n}$ for $n\ge10$ and
$m\le20$, so we may apply Theorem \ref{thetheorem}(b) with
$\mu=21$ and $n\ge 16$. (Here the only $k\ne0$ allowed
is for $k=m=1$, but it is known and easy to see that 
a double point and general points of multiplicity 1
impose independent conditions on forms on \pr2 of degree $\alpha$.
Thus $(\alpha+3/2)^2/2>{\alpha+2 \choose 2}>3+(n-1)=n+2$, so
for $k=m=1$,
$\alpha_0((m^{[n-1]},m+k)) \ge {{mn+k}\over{n}} \sqrt{(n-1/\mu)}$
since $(\alpha_0+3/2)^2\ge(\alpha+3/2)^2>2n+4\ge(\sqrt{n}+2)^2$
for $n\ge16$, and $(\sqrt{n}+2)^2\ge ({n+1\over \sqrt{n}}+3/2)^2
=({{mn+k}\over{n}} \sqrt{n}+ 3/2)^2>({{mn+k}\over{n}} \sqrt{(n-1/\mu)}+3/2)^2$.)
We thus immediately obtain an explicit bound which for most
$n$ is substantially better than what was known
previously:\footnote{After submission of this paper, a result for
$m\le 42$ has been announced by M. Dumnicki \cite{refDu} which,
together with \cite{refCM} for the $k\ne 0$ case, imply the stronger
bound $\sc{n}\ge(\sqrt{1/n}\,)\sqrt{1-1/43\,n}$ if $n\ge 16$.}

\begin{cor}\label{corRem}
For every $n\ge 16$, 
$$\sc{n}\ge \hbox{max}\Bigg({1\over \sqrt{n}}\sqrt{1-{1\over n(1+\lfloor\sqrt{n}\rfloor(\lfloor \sqrt{n}\rfloor-3)/2)}},\,
{1\over \sqrt{n}}\sqrt{1-{1 \over{21\,n}}}\,\,\Bigg).$$
\end{cor}

As a final application, again for blow ups $Y$ of $X=\pr2$ where $L$ is a line, we  
obtain an improved criterion for which divisor classes of the form
$dL-m(E_1+\cdots+E_n)$ are ample. If Nagata's conjecture \cite{refN}
is true, it is not hard to see that $F=dL-m(E_1+\cdots+E_n)$ is ample whenever
$d$ and $m$ are positive integers such that $d^2>m^2n$, where
$\pi: Y\to \pr2$ is given by blowing up $n\ge10$ very general points
and $L$ is the class of a line.
That $F$ is in fact ample has been verified for $m=1$ \cite{refXc}, $m=2$ \cite{refBi} and 
$m=3$ \cite{refT}. Our result extends these substantially for large $n$
(see \cite{refH}, however, for an even stronger result if one merely wishes to conclude
that $F$ is nef):

\begin{cor}\label{simplifiedampCor}
Let $n\ge 16$, $t>\sqrt{n}m$, and $m>0$
be integers and
consider the divisor class $F=tL-m(E_1+\cdots+E_n)$ on the blow 
up $Y$ of \pr2 at $n$ general points, where $L$ is
the pullback to $Y$ of a line in \pr2. 
If $1\le m<\sqrt{\lfloor\sqrt{n}\rfloor(\lfloor \sqrt{n}\rfloor-3)/2 + 1-1/n}$, then
$F$ is ample.
\end{cor}

We end this introduction by discussing
Corollary \ref{corRem} in the context of what was known previously in case $X=\pr2$.
It is convenient for comparison to express lower bounds
for Seshadri constants on \pr2 in the form  $(\sqrt{1/n}\,)\sqrt{1-1/f(n)}$.
Note that the larger $f(n)$ is, the better is the bound.
Perhaps the best previous general bound 
is given in \cite{refT}, for which $f(n)=12n+1$ for all $n\ge10$. 
For Corollary \ref{corRem}, which applies for  all $n\ge17$, $f(n)$ can be taken to be
quadratic in $n$ but always larger than $12n+1$.

The article \cite{refBi} gives bounds which for special values of $n$ are
better than those of \cite{refT}, and for these special values $f(n)$ is
quadratic in $n$. (In particular, if
$n=(ai)^2\pm 2i$ for positive integers $a$ and $i$,
then $f(n) = (a^2i\pm1)^2$, and, if $n=(ai)^2+i$ for 
positive integers $a$
and $i$ with $ai\ge 3$, then $f(n) = (2a^2i+1)^2$).)
However, except in special cases, such as 
when $n-1$ or $n\pm2$ is a square, the bounds of
Corollary \ref{corRem} are better for $n$ large enough.  (To see this 
look at coefficients of the $n^2$ term in $f(n)$.)

Bounds are also given in \cite{refH}; they apply for all values of $n$
for all surfaces and are almost always better than any bound for which $f(n)$ is linear 
in $n$ (more precisely, given any constant $a$, let $\nu_a(n)$ be 
the number of integers
$i$ from 1 to $n$ for which $f(i)$ from \cite{refH} is bigger than
$ai$; then $\hbox{lim}_{n\to\infty}\nu_a(n)/n=1$). 
However, although the bounds in \cite{refH} are not hard to 
compute for any given value of $n$, there is no simple explicit formula for 
$f(n)$, so it is hard to make general comparisons.
Nonetheless, computations in case $X=\pr2$ for specific values of $n$ 
suggest that the bounds we obtain here for \pr2
are typically if not almost always better than those of \cite{refH}.

It is worth noting that the bounds in Corollary \ref{corRem} are not the best that
one can obtain using our results here in conjunction with the methods of \cite{refHR}.
While \cite{refHR} does give explicit formulas that hold in general,  
applying the methods of \cite{refHR} for specific values of $n$ usually gives notably 
better results than one can express in terms of an explicit formula.
Since the simple explicit formula for $f(n)$ as given in Corollary \ref{corRem}
is based on an explicit but necessarily suboptimal formula from \cite{refHR},
one can usually get better results for specific values of $n$
by directly applying the methods of Section \ref{submax} and 
\cite{refHR}. (For specific examples of this, see the unpublished posting
\cite{refoldpreprint}.)

\section{Main Results}\label{submax}

In the first section we obtain results about abnormal
curves in general. 
In the second section we sharpen and apply those results
in the homogeneous case. For the rest of this paper we assume
that $X$ is a smooth projective surface. 
 
\subsection{Abnormal Curves}\label{disc}

Let $\pi:Y\to X$ be
obtained by blowing up distinct points $p_i$ on $X$ and let $E_i=\pi^{-1}(p_i)$.
Let $L$ be a nef divisor on $X$.
Abnormality, as we introduced it above, is related to
nefness of divisors on $Y$ of the form $\pi^*L-E_1-\cdots-E_n$.
In order more generally to study nefness of
divisors of the form $\pi^*L-l_1E_1-\cdots-l_nE_n$,
it is convenient to extend our notion of abnormality. 
Let $F$ be a numerical equivalence divisor class on $Y$.
We will then say a curve $D\subset Y$ is 
$F$-{\it abnormal\/} if $D$ is reduced and irreducible 
with $F\cdot D<0$. In case the points $p_i$ are general, $L$ is nef on $X$ and
$F=\pi^*L-(E_1+\cdots+E_n)$, then a curve $C\subset X$ is $L$-abnormal
according to our previous use of the word, if and only if
its proper transform $\tilde C$ is $F$-abnormal.

For simplicity, we will by identification just write $L$ in place of $\pi^*L$.
The next lemma establishes a connection between values of $s$ for
which $F_s=L-s(l_1E_1+\cdots +l_nE_n)$ is nef and the occurrence of
abnormal curves.

\begin{lem}\label{introlem}
Let $L$  be a nef divisor on
$X$, let $\pi:Y\to X$ be obtained by blowing up $n$ distinct points $p_i$ on $X$
and let $F_t=L-t(l_1E_1+\cdots+l_nE_n)$, where $E_i=\pi^{-1}(p_i)$
and where $t$ and each $l_i\ge0$ is real
(such that ${\bf \ell}= (l_1,\ldots,l_n)$ is not 0). 
\begin{itemize}
\item[(a)] If $F_t$ is nef, then $0\le t\le \sqrt{L^2/{\bf \ell}^2}$.
\item[(b)] Let $0\le t\le\sqrt{L^2/{\bf \ell}^2}$. If $D$ is an
$F_t$-abnormal curve on $Y$, then
the largest $s$ such that $F_s$ is nef is at most
$(L\cdot D)/D\cdot(l_1E_1+\cdots+l_nE_n)$. Moreover,
any such $D$ satisfies $D^2<0$.
\item[(c)] Let $t$ and $\ell$ be rational and $0\le
t<\sqrt{L^2/{\bf \ell}^2}$. Then  
the following are equivalent:
\begin{itemize}
\item[(i)] there exists a numerical equivalence class $H$ which
for general points $p_i$ is the class of an $F_t$-abnormal curve;
\item[(ii)] $\sc{X,L,n,{\bf \ell}}<t$; and
\item[(iii)] $F_t$ is not nef for any choice of the points $p_i$.
\end{itemize}
\end{itemize}
\end{lem}

\begin{proof} (a) We have $0\le t$ since $F_t$ is nef and hence $tl_i=F_t\cdot E_i\ge0$ for all $i$, while
$t\le\sqrt{L^2/{\bf \ell}^2}$ follows since any nef divisor has non-negative self-intersection.

(b) If $F_t$ is not nef, then $L^2>0$ (else $t=0$ and $F_t=\pi^*L$ is nef).
Since $F_t$ is not nef, there is an $F_t$-{\it abnormal\/} curve $D$. 
If $F_s$ is nef, then $L\cdot D-s(l_1E_1+\cdots+l_nE_n)\cdot D=
F_s\cdot D\ge0$, so $s\le (L\cdot D)/D\cdot(l_1E_1+\cdots+l_nE_n)$.

To see $D^2<0$, note that up to numerical equivalence,
we can write $D$ as $C'-m_1E_1-\cdots-m_nE_n$, for some integers
$m_i$ where $C'=\pi^{-1}(\pi(D))$. Since $t\ge0$, we have 
$F_t\cdot E_i\ge0$ for all $i$, so $D$ cannot by $E_i$ for any $i$.
Thus $\pi(D)$ is a curve, and $m_i\ge0$ for each $i$.
Since $L^2>0$, we can by the Hodge index theorem write
$C=cL+B$ for some real $c\ge0$ and some ${\bf R}$-divisor
$B$ with $B\cdot L=0$ and $B^2\le0$, where $C=\pi(D)$.
Thus $C^2=c^2L^2+B^2\le (cL)^2 = (C\cdot L)^2/L^2 < (l_1m_1+\cdots+l_nm_n)^2/{\bf \ell}^2$,
where the strict inequality follows since $D$ is $F_t$-abnormal and $t\le \sqrt{L^2/{\bf \ell}^2}$.
But $(l_1m_1+\cdots+l_nm_n)^2/{\bf \ell}^2\le \sum_im_i^2$ by Cauchy-Schwarz,
so $D^2=C^2-\sum_im_i^2<0$, as claimed.

(c) If an $F_t$-abnormal curve of class $H$ exists for general
sets of distinct points $p_i$, then since 
$\sc{X,L,n,{\bf \ell}}=\sc{X,L,n,l_1p_1, \dots, l_np_n}$ on a dense set,
from the definitions it follows that  
$\sc{X,L,n,{\bf \ell}}\le L\cdot H/(H\cdot (l_1E_1+\cdots+ l_nE_n)) < t$.
If $\sc{X,L,n,{\bf \ell}}<t$, then by definition
$F_t$ is not nef for every set of points $p_i$.
Finally, if $F_t$ is not nef for every set of points ${\bf
p}=(p_1,\ldots,p_n)$, then  
for each choice of the points ${\bf p}$ one can choose an $F_t$-abnormal
$H_{\bf p}$. By Lemma \ref{genblowups}, there are
only finitely many 
classes of such $H_{\bf p}$ in $\hbox{Pic}(Y)/\hbox{$\sim$}$, and each
of them is effective on a Zariski-closed set. Hence one of them (say
$H$) must be effective for all choices of the points ${\bf p}$ and
irreducible for a general set of 
points ${\bf p}$, with $H\cdot F_t<0$. \end{proof}

We now state a lemma of particular interest, since it applies
to the case of $n$ general points on $X=\pr2$.

\begin{lem}\label{introlemb}
Assume the hypotheses of Lemma \ref{introlem}
together with the additional hypothesis that the points $p_i$ are general points of $X$.
If $D$ and $F_t$ are as in Lemma \ref{introlem}(b) with $1=l_1=\cdots=l_n$, and if 
every ${\bf R}$-divisor on $X$ (up to
numerical equivalence) is a real multiple of $L$, then 
the largest $s$ such that $F_s$ is nef is precisely 
$s = (L\cdot D)/D\cdot(E_1+\cdots+E_n)$; i.e.,
$\sc{X,L,n} = s$.
\end{lem}

\begin{proof} We use the argument of Proposition 4.5 of \cite{refS}.
Suppose that there is another $F_t$-abnormal curve $D'$,
whose class is $C''-m'_1E_1-\cdots-m'_nE_n$.
Since the points $p_i$ are general, we may assume
that $m_1\ge m_2\ge\cdots\ge m_n$ and $m_1'\ge m_2'\ge \cdots\ge m'_n$, and 
so by Chebyshev's sum inequality we have
$((m_1+\cdots+m_n)/n)((m_1'+\cdots+m_n')/n)\le (m_1m_1'+\cdots+m_nm_n')/n$.
But $C$ and $C'$ are positive multiples of $L$, so there are positive reals
$c$ and $c'$ such that $C=cL$ and $C'=c'L$. 
We have therefore that $cL^2/(m_1+\cdots+m_n)$ and $c'L^2/(m'_1+\cdots+m'_n)$
both are less than $\sqrt{L^2}/\sqrt{n}$, and hence
$${ncc'(L^2)^2\over \sum_im_im'_i} \le {cc'(L^2)^2\over {\sum_im_i\over n}{\sum_im'_i\over n}}
< {n^2L^2\over n} = nL^2,$$
so $cc'(L^2) < \sum_im_im'_i$; i.e., $D\cdot D'<0$.
Since $D$ and $D'$ are integral, we must have $D=D'$.
Thus every $F_t$-abnormal curve $B$ gives the same value
for $(L\cdot B)/B\cdot(E_1+\cdots+E_n)$. 
By (b), $F_s$ cannot be nef for any value of $s$ bigger than 
$s=(L\cdot D)/D\cdot(E_1+\cdots+E_n)$,
yet for this value of $s$ we have just shown there are no $F_s$-abnormal curves,
so $F_s$ is in fact nef, and hence $\sc{X,L,n} = s$.
\end{proof}

\eatit{Brian, I have not had time to look into this.
A nonhomogeneous version of part (d) doesn't seem to work.
Consider F=3.61L-3E_1-2E_2,
D=41L-30E_1-30E_2,
D'=3L-3E_1-E_2, where L^2=1.
Then FD<0, FD'<0, F^2>0, but DD'>0.

However, a nonhomogeneous version of Chebyshev is true:
if a_1 >= a_2 >= ... >= a_n >= 0,
b_1 >= b_2 >= ... >= b_n >= 0, and
l_1 >= l_2 >= ... >= l_n > 0,
let a = (a_1l_1+...+a_nl_n)/(l_1+...+l_n) be the weighted average a of the a's,
b the weighted average of the b's, and c the weighted average of the
products a_ib_i. Then ab <= c.

Proof: We may assume (l_1+...+l_n)=1.
Then ab = (ba_1l_1+...+ba_nl_n) <= b_1a_1l_1 + b*_1(a_2l_2+...+a_nl_n),
where we define b*_1 by
b_1l_1+b*_1(l_2+...+l_n)=b=(b_1l_1+...+b_nl_n). The inequality is because
we are increasing bl_1 to turn it into b_1l_1 by the same amount we are decreasing b(l_2+...+l_n)
to turn it into b*_1(l_2+...+l_n) [i.e., (b_1-b)l_1=(b-b*_1)(l_2+...+l_n)], 
but a_1 >= a_2 so (b_1-b)l_1a_1 >= (b-b*_1)(l_2+...+l_n)a_2 >= (b-b*_1)(l_2a_2+...+l_na_n).
But b*_1=(b_2l_2+...+b_nl_n)/(l_2+...+l_n). Now induct.
}

To state the general fact used in
Lemma \ref{introlem}(c), we define the notion of a {\it sufficient test system}.
Let ${\bf p}=\{p_1,\ldots,p_n\}$ be a set of distinct points on a surface $X$,
and let $\pi:Y_{\bf p}\rightarrow X$ be the morphism obtained by blowing 
up the points $p_i$ with, as usual, $E_i=\pi^{-1}(p_i)$. 
Given a $\bf Q$-divisor $L$ on $X$ and nonnegative rationals 
$m_1,\ldots,m_n$, consider a set $\{D_1,\ldots, D_k\}$ of 
numerical equivalence classes of divisors on $X$ 
together with vectors ${\bf h}_1,\dots, {\bf h}_k \in {\bf Z}^{n}_{\ge 0}$.
We refer to $\{(D_i,{\bf h}_i)\}_{i=1,\dots,k}$
as an $(L, \{m_i\})$-sufficient test system
if whenever ${\bf p}$ is such that 
none of the classes $C_i=D_i - h_{i1} E_1 - \cdots -h_{in} E_n$ is
(up to numerical equivalence) the class
of a reduced irreducible curve, then $F=L - m_1E_1-\cdots-m_nE_n$
is nef. Remark that by definition a $(L, \{m_i\})$-sufficient test
system is always finite.

\begin{lem}\label{genblowups}
Let $L$ be a big and nef  
$\bf Q$-divisor on $X$, and let $m_1, \dots, m_n$ be non-negative rationals with
$m_1^2+\cdots+m_n^2<L^2$. Then there exists an
$(L, \{m_i\})$-sufficient test system $\{(D_i,{\bf h}_i)\}_{i=1,\dots,k}$.
Moreover, if $U\subset X^n$ is the set of all $n$-tuples
of distinct points, then for each class 
$C_i=D_i - h_{i1} E_1 - \cdots -h_{in} E_n$ the subset of $U$ such that
$C_i$ is the class of an effective divisor on the blowup
of ${\bf p}\in U$ is Zariski-closed. (In particular, the subset of 
$U$ such that $F$ is nef on the blowup
of ${\bf p}\in U$ is Zariski-open.)
\end{lem}

\begin{proof} Clearly, there is
an $s$ such that $sF$ is effective.
Let $L_1, \ldots, L_\rho$ be ample effective divisors which generate the
group of numerical equivalence classes on $X$. For suitable
$a_{i0}, a_{ij}\in {\bf N}$, with $1\le i\le \rho, 1\le j\le n$, the divisor
classes $A_{i0}=a_{i0}L_i-(E_1+\cdots+E_n)$, $A_{ij}=
a_{ij}L_i-(E_1+\cdots+E_n)-E_j$ are ample and effective, and they
generate the group of numerical equivalence classes, independently of the
choice of the points. Let $d_{ij}=sF\cdot A_{ij}$ for all $i$ and all $0\le j\le n$. If $C$ is a
divisor such that both $|C|$ and $|sF-C|$ are nonempty 
(which is necessary in order to have an $F$-abnormal curve $C$),
then $0< C\cdot A_{ij}\le d_{ij}$. Moreover, the class of an
irreducible curve meeting $F$ negatively must be of the form 
$C=D - h_{1} E_1 - \cdots -h_{n} E_n$, and clearly there are only a
finite number of numerical equivalence classes of such $C$
satisfying $0< C\cdot A_{ij}\le d_{ij}$. Let these classes be 
$C_i=D_i - h_{i1} E_1 - \cdots -h_{in} E_n$, $i=1, \dots, k$; we have
shown that $\{(D_i,{\bf h}_i)\}_{i=1,\dots,k}$ is
a $(L, \{m_i\})$-sufficient test system.

The set $\{{\cal H}_\gamma\}_{\gamma \in \Gamma_i}$ of all
components in the Hilbert scheme of curves in $X$ 
numerically equivalent to $D_i$  is indexed by some finite set
$\Gamma_i$ (see e.g. \cite{refMu}, lecture 15).
Since there are only finitely many $D_i$, it follows that $\Gamma=\bigcup \Gamma_i$
is finite. For each $\gamma \in \Gamma$, there is a flat family
$\phi_\gamma:{\cal D}_\gamma\subset X\times {\cal H}_\gamma\rightarrow {\cal
H}_\gamma$ whose members are the curves parameterized by ${\cal
H}_\gamma$; every $F$-abnormal curve of class $C_i$ 
occurs as the birational transform of a fiber of some $\phi_\gamma$,
$\gamma \in \Gamma_i$, 
which has multiplicity $h_{ij}$ at a point $p_j \in X$. Now the sets of
(distinct) points $(p_1, \dots, p_j)\in U$ such that there exists a fiber of
$\phi_\gamma$ with multiplicity at least $h_{ij}$ at the point $p_j$ is
Zariski closed in $U$ (an explicit construction of
this closed set, using sheaves of principal parts, can be found
e.g. in \cite{refKP}, section 4). Since the subset of $U$ such that
$C_i$ is the class of an effective divisor on the blowup
of ${\bf p}\in U$ is the union of the finitely many closed subsets
determined by the $\phi_\gamma$, $\gamma \in \Gamma_i$, it follows
that it is Zariski-closed.

 Finally, the divisor
$F$ is nef if and only if
none of the classes $C_i$ is effective,
and we have seen that the set of points $p_i$ for which
none of them is effective is open.
\end{proof} 

Such general claims as in Lemma \ref{genblowups}
regarding the existence of a finite set of test classes
for $F_t$ to be nef can be sharpened and made more explicit
in the case of general blow-ups, as we now show.

Given a big and nef divisor $L\subset X$ and non-negative integers
${\bf \ell}=(l_1, \dots, l_n)$, let $F=dL-l_1E_1-\cdots-l_nE_n$
where $d=\sqrt{{\bf \ell}^2/L^2}$, so $F^2=0$.
For each real ${\delta}\ge0$, consider the ${\bf R}$-divisor 
$F(\delta)=d'L-l_1E_1-\cdots-l_nE_n$
where $d'=\sqrt{({\bf \ell}^2+\delta)/L^2}$; note that
$F(\delta)^2=\delta$.
 The next lemma can be seen as a sharpening and extension of
theorem 4.1 in \cite{refBa} to the case of multipoint Seshadri constants:

\begin{lem}\label{testLem}
Let 
$\pi: Y\to X$ be the blow up of general points $p_1,\ldots,p_n\in X$.
Let $F$ and $F(\delta)$ be as in the preceding paragraph
with $\delta > 0$. 
If $H$ is the class of an $F(\delta)$-abnormal curve $\tilde C$, 
then $H=\pi^* C-h_1E_1-\cdots-h_nE_n$
for some non-negative integers $h_1,\ldots,h_n$ and for some effective divisor
class $C$ on $X$ such that:
\item{(a)} $h_1^2+\cdots+h_n^2<(1+d^2L^2/\delta)^2/\gamma$, 
where $\gamma$ is the number of nonzero coefficients $h_1,\ldots,h_n$, and
\item{(b)} $h_1^2+\cdots+h_n^2-a\le C^2 \le (C\cdot L)^2/L^2 < 
(l_1h_1+\cdots+l_nh_n)^2/(d^2L^2+\delta)$,
where $a$ is the minimum positive element of $\{h_1,\ldots,h_n\}$.
\end{lem}

\begin{proof} 
The class $H$ of $\tilde C$ must 
be of the form $H=\pi^* C-h_1E_1-\cdots-h_nE_n$, with $C$ effective (since
$\tilde C$ is effective) and each $h_i$ non-negative 
(since $\tilde C$ is irreducible
and $F(\delta)\cdot E_i\ge0$ holds for all $i$).

First consider (b). By \cite{refXb}, Lemma 1, we have $\tilde
C^2\ge-a+1$ if $a>1$. It is easy to see that $\tilde C^2\ge-1$ if $a=1$, 
for suppose $\tilde C\cdot E_i=1$ yet $\tilde C^2<-1$. 
Then we would have $(\tilde C + E_i)^2<0$, hence
$|\tilde C + E_i|$ is fixed.
However, the linear system $|\tilde C + E_i|$ corresponds to a complete
linear system on the surface $Y'$ obtained by contracting $E_i$;
$|\tilde C|$ corresponds to the subsystem vanishing at $p_i$.
Since $p_i$ is a general point, $|\tilde C + E_i|$ cannot be fixed,
which contradicts $\tilde C^2<-1$ when $a=1$.
Hence we may assume $\tilde C^2\ge -a$,
so $h_1^2+\cdots+h_n^2-a\le C^2$. Also, since $L$ is
big and nef, the index theorem (as in the proof of Lemma \ref{introlem}(b))
gives  $C^2 L^2\le (C \cdot L)^2$.
On the other hand, $F(\delta)\cdot
\tilde C<0$ gives $(C\cdot L)^2 < (l_1h_1+\cdots+l_nh_n)^2L^2/(d^2L^2+\delta)$.

Now consider (a). Let $h=\sqrt{h_1^2+\cdots+h_n^2}$.
From (b) we have
$h^2-a < (l_1h_1+\cdots+l_nh_n)^2/(L^2d^2+\delta)\le 
d^2L^2h^2/(d^2L^2+\delta)$, so $h^2 < d^2L^2h^2/(d^2L^2+\delta) + a$. 
But $a^2\le h^2/\gamma$, so
we have $h^2 < d^2L^2h^2/(d^2L^2+\delta) + h/\sqrt{\gamma}$, and solving 
for $h$ gives the result.
\end{proof}

For each $\delta>0$, let $O_n(F(\delta))$ be the set of all numerical 
equivalence classes of divisors $H=\pi^*C-h_1E_1-\cdots-h_nE_n$ 
where $C$ is the class of an effective divisor on $X$ and $C$ and the $h_i$
satisfy the inequalities in
Lemma \ref{testLem}(a, b). Then $O_n(F(\delta))$ is the set of obstructions to 
$F(\delta)=d'L-(l_1E_1+\cdots+l_nE_n)$ being nef;
i.e., $O_n(F(\delta))$  contains the class of every $F(\delta)$-abnormal curve (if any). In particular, $O_n(F(\delta))$ is an $(L,\{l_1,\ldots,l_n\})$-sufficient test system.
Let $o_n(F(\delta))$ be the set of ratios $L\cdot C/(l_1h_1+\cdots+l_nh_n)$ for
all $H\in O_n(F(\delta))$. 

\begin{thm}\label{maindiscthm}
Let $L$, $F(\delta)$, $Y$ and $X$ be as in
Lemma \ref{testLem}. Then $o_n(F(\delta))$ is a
finite set for each $\delta>0$, and the union 
$U_n=\cup_{\delta>0}o_n(F(\delta))$ is discrete,
with $t=\sqrt{L^2/{\bf \ell}^2}$ as the unique limit point (if any).
Moreover, if $F(\delta)$ is not nef for some $\delta>0$
(which is equivalent to $\varepsilon(X,L,n,{\bf \ell})<\sqrt{L^2/{\bf \ell}^2}$),
then $\varepsilon(X,L,n,{\bf \ell})$ is the maximum $t$ such that
$F_t=L-t(l_1E_1+\cdots+l_nE_n)$ is nef and this $t$ is an element of $o_n(F(\delta))$;
i.e., $\varepsilon(X,L,n,{\bf \ell})\in U_n$.
\end{thm}

\begin{proof} Lemma \ref{testLem} implies that $o_n(F(\delta))$ is finite.
If $\delta'<\delta$, then every element $t$ of $o_n(F(\delta'))$ not in
$o_n(F(\delta))$ is bigger than every element of $o_n(F(\delta))$; in particular,
$\sqrt{L^2/({\bf \ell}^2+\delta)}\le t < \sqrt{L^2/({\bf \ell}^2+\delta')}$, hence
the only possible limit point is $t=\sqrt{L^2/{\bf \ell}^2}$.
Note that $(1/c)F_c=F(\delta)$ exactly when 
$\delta=L^2/c^2-{\bf \ell}^2$, so if $\delta=L^2/c^2-{\bf \ell}^2$, 
then $F(\delta)$ is nef if and only if $F_c$ is, so
$F(\delta)$ not being nef for some $\delta>0$
is by Lemma \ref{introlem}(c) equivalent to $\varepsilon(X,L,n,{\bf \ell})<\sqrt{L^2/{\bf \ell}^2}$.
If $F(\delta)$ is not nef, take $t$ to be the infimum for
$L\cdot C/(l_1h_1+\cdots+l_nh_n)$ over all 
classes $H=C-(h_1E_1+\cdots+h_nE_n)$
of $F(\delta)$-abnormal curves. Thus $t\in o_n(F(\delta))$
since $o_n(F(\delta))$ is finite, and 
$L-t(l_1E_1+\cdots+l_nE_n)$ is nef since
we have chosen $t$ small enough to eliminate
all obstruction classes. Finally, by Lemma \ref{introlem}(c), 
we also have $\varepsilon(X,L,n,{\bf \ell})=t$.
\end{proof}

\eatit{
If $L$ is not ample, then $O_n(F(\delta))$ is infinite.
Even if $L$ is ample, $O_n(F(\delta))$ can be infinite,
since there can be some non-effective $B$ with $L\cdot B=0$
and hence $C+\lambda B-(h_1E_1+\cdots+h_nE_n)\in O_n(F(\delta))$
for all $\lambda$. If $L$ were ample and if
we had a set of ample generators $A_1,\ldots,A_\rho$
for $\hbox{Pic}(Y)/\hbox{$\sim$}$, we could cut 
$O_n(F(\delta))$ down to a finite set by only taking elements
$H\in O_n(F(\delta))$ satisfying $sF(\delta')\cdot A_i \ge H\cdot A_i>0$
where $0<\delta'\le\delta$ is chosen so $F(\delta')$ is rational
and $s$ is chosen so that $sF(\delta')$ is integral and effective.
}

Observe that from Lemma \ref{introlem}(c) and (d) it follows that 
$\sc{X,L,n,p_1, \dots,p_n}=\sc{X,L,n}$ for general points whenever
$\sc{X,L,n}<\sqrt{L^2/n}$ and the group of numerical
equivalence classes has rank one. 
However, by Theorem \ref{maindiscthm} it now follows for all
${\bf \ell}$ and all $X$ that
$\sc{X,L,n,{\bf\ell}}=\sc{X,L,n,l_1p_1, \dots, l_np_n}$ for 
general points whenever $\sc{X,L,n,{\bf \ell}}<\sqrt{L^2/\ell^2}$.
To see this, let $t=\sc{X,L,n,{\bf \ell}}$.
By Lemma \ref{introlem}(c), $F_t$ is nef for some choice of points
$p_i$, and hence by Lemma \ref{genblowups} for an open
set. Thus on some nonempty open set we have
$\sc{X,L,n,l_1p_1, \dots, l_np_n}\ge t$.
On the other hand, by the discreteness claim of 
Theorem \ref{maindiscthm} there exists a $t'$ such that
$t'>t$ but such that no element of $\cup_\delta o_n(F(\delta))$
is in the interval $(t,t']$. By Lemma \ref{introlem}(c)
it follows that there is an open set for which
there exists an $F_{t'}$-abnormal $H$. 
Since $F_t\cdot H\ge0$ but $F_{t'}\cdot H<0$,
it must be that $H\cdot L/(H\cdot (l_1E_1+\cdots+l_nE_n))$
is in the interval $[t,t')$, and hence that
$t=H\cdot L/(H\cdot (l_1E_1+\cdots+l_nE_n))$.
Thus on this nonempty open set we also have
$t\ge \sc{X,L,n,l_1p_1, \dots, l_np_n}$.

\subsection{Applications}\label{appl}

We now turn our attention to obtaining explicit bounds
on homogeneous Seshadri constants. 
We begin this section by describing our conceptual basis 
for bounding Seshadri constants. Given general
points $p_i\in X$, $1\le i\le n$ on $X$ and 
a big and nef divisor $L$ on $X$,
let $\pi:Y\to X$ be obtained from $X$ by blowing up the points.
Then $\sc{X,L,n}\ge t$ whenever $F_t=L-t(E_1+\cdots+E_n)$ is 
big and nef, by Lemma \ref{introlem}(c)
(the case that $t$ is real follows by taking the limit of smaller rational values). 

In order to show $F_t$ is nef for a given $t$ for which
$F_t^2>0$, we first consider the set $O_n(F_t)$ of test classes,
which we obtained from Lemma \ref{testLem}.
We can explicitly determine the finite set $o_n(F_t)$.
If each test class is shown not to
be the class of a reduced, irreducible curve (by showing, for example, 
that none is the class of an effective divisor), it follows that $F_t$ is nef
and hence that $\sc{X,L,n}\ge t$.
However, Lemma \ref{testLem} applies more generally 
to classes $F=L-t(l_1E_1+\cdots+l_nE_n)$.
Since hereafter we will focus on $F=L-t(E_1+\cdots+E_n)$,
it behooves us to make better use of the fact
that the coefficients $l_i$ are equal. Doing so allows us to
significantly sharpen Lemma \ref{testLem}, which we state as
Corollary \ref{SzCor}.

We need the following lemma, which generalizes a result of \cite{refS}:

\begin{lem}\label{LemSz}
Let $F$ be an ${\bf R}$-divisor class on $X$ with $F\cdot L>0$ for some
big and nef class $L$ and with $F^2\ge 0$. Let 
$C_1,\ldots,C_r$ be distinct $F$-abnormal 
curves. Then up to numerical equivalence 
their divisor classes $[C_1],\ldots,[C_r]$ are linearly
independent in the divisor class group on $X$.
\end{lem}

\begin{proof} If $[C_1],\ldots,[C_r]$ are 
dependent, we can find a nontrivial non-negative integer combination 
$D$ of some of the classes $[C_1],\ldots,[C_r]$
and another nontrivial non-negative integer combination $D'$ of the rest 
of the classes $[C_1],\ldots,[C_r]$, such that, up to numerical equivalence, $D=D'$. 
But $F\cdot D<0$, so for some real number $\delta>0$
we must have $(F+\delta L)\cdot D=0$ with $(F+\delta L)^2>0$, hence by the index theorem
we must have $D^2<0$, which contradicts $D^2=D\cdot D'\ge0$.\end{proof}

The analysis of what $F$-abnormal curves can occur is especially simple
when the coefficients $F\cdot E_i$ are all equal. In particular, as our 
next result generalizing and extending methods and results of 
\cite{refX}, \cite{refSS} and \cite{refR} shows,
they must be almost uniform, where we call a class of the 
form $\pi^* C-m(E_1+\cdots+E_n)$ {\it uniform}, and
we call a class of the form
$\pi^* C-m(E_1+\cdots+E_n)-kE_i$ {\it almost uniform}
(called almost homogeneous in \cite{refSS}).

\begin{cor}\label{SzCor}
Let $L$ be a big and nef divisor on $X$. 
Let $\pi:Y\rightarrow X$ be the blow
up of $n\ge1$ general points $p_1,\ldots,p_n\in X$. Consider  the ${\bf
R}$-divisor class $F=(\sqrt{n/L^2})\pi^*L-E_1-\cdots-E_n$, and let $H$
be a divisor class on $Y$ with $F\cdot H<0$. If $H$ is the class 
of an $F$-abnormal curve,
then there are integers 
$m>0$, $k$ (where we require $k=0$ if $n=1$)
and $1\le i\le n$ and an effective divisor $C$ on $X$ such that:
\begin{itemize} 
\item[(a)] $H=\pi^*C-m(E_1+\cdots+E_n)-kE_i$;
\item[(b)] either $k>-m$ and $k^2< (n/(n-1))\,\hbox{min}\,(m,
m+k)$, or $m=-k=1$;
\item[(c)] $(m^2n+2mk+\hbox{max}(k^2-m, k^2-(m+k), 0))L^2
\le C^2L^2 \le (C\cdot L)^2<(m^2n+2mk+k^2/n)L^2$ when $k^2>0$, but
$(m^2n-m)L^2\le C^2 L^2 \le (C\cdot L)^2<(m^2n)L^2$ when $k=0$; and
\item[(d)] $C\cdot(C+K_X) - (m+k)^2 - (n-1)m^2 + mn + k\ge -2$.
\end{itemize}
\end{cor}

\begin{proof} The case $n=1$ (and so $k=0$) is easy to treat along the same lines as below;
we leave it to the reader. Thus we assume $n\ge2$. 

(a) In \cite{refSS}, corollary 2.8,  this result is proved for
surfaces of Picard number 1. We adjust their argument to prove the
result for arbitrary Picard numbers. 
Because the points are general and $F$ is uniform, permuting the 
coefficients $m_i$ of the class $H=\pi^*C-m_1E_1+\cdots+m_nE_n$ of an
$F$-abnormal curve gives another such class.  
Since all such permutations are in the subspace of the span of
$\pi^*C,E_1,\dots,E_n$ orthogonal 
to $F-(F\cdot H)/(C\cdot L)\pi^*L$, it follows from Lemma \ref{LemSz} that there 
are at most $n$ such curves. But it is not hard to check that
there are always more than $n$ permutations unless at most
one of the coefficients is different from the rest.
Thus $H$ is of the form $H=\pi^*C-m(E_1+\cdots+E_n)-kE_i$
with $1\le i\le n$, which gives (a).

Since $H\cdot F(0)=H\cdot F<0$, it follows that for $\delta>0$ small
enough, $H\cdot F(\delta)<0$. For the proof of (b) and (c), fix a
$\delta>0$ such that $H$ is the class of a $F(\delta)$-abnormal curve.

Consider (b). Since $H$ is the class of a reduced irreducible curve
with $C\cdot L>0$, we must have $H\cdot E_i\ge 0$ for all $i$,
hence $-m\le k$. If $k=-m$, then Lemma \ref{testLem}(b) says
$(m^2(n-1)-m)<(m(n-1))^2/n$, which simplifies to $m^2(n-1)<mn$, and hence
$m=-k=1$. 
Now, again by Lemma \ref{testLem}(b) with $a=\hbox{min}\,(m, m+k)$,
we have $(m^2n+2mk+k^2 - a)n < (mn+k)^2$,
which simplifies to give $k^2< (n/(n-1))\,(a)$.

Likewise,
(c) follows from Lemma \ref{testLem}(b) in the case that $k=0$, as does
$(m^2n+2mk+\hbox{max}(k^2-m, k^2-(m+k)))L^2\le (C\cdot
L)^2<(m^2n+2mk+k^2/n)L^2$ when  
$k\ne0$. If $k\ne0$, then $\pi^*C-m(E_1+\cdots+E_n)-kE_1$ 
and $\pi^*C-m(E_1+\cdots+E_n)-kE_n$
are classes of distinct irreducible curves, so their intersection
is non-negative, hence $m^2n+2mk\le C^2$, and
$(m^2n+2mk+\hbox{max}\,(k^2-m,k^2-(m+k),0))L^2\le (C\cdot L)^2$ as claimed.

Finally, we prove (d). A reduced, irreducible curve 
must have a non-negative genus $g$,
hence by adjunction we must have $H^2+K_Y\cdot H = 2g(H)-2\ge -2$, 
which is (d). \end{proof}
 
It may be interesting to note that item (d) above is implied by (b)
and (c) if $X=\pr2$ and the number of 
points is $n\ge 11$. The proof of 
this implication follows from a straightforward but somewhat lengthy
computation that we leave to the interested reader to carry through.

It may also be of interest that 
Corollary \ref{SzCor} takes the following very simple form
if $m<n$. Since we will not use the following result we omit a proof.

\begin{cor}\label{lemAAA}
Let $\pi:Y\rightarrow X$ be the blow
up of $n$ general points $p_1,\ldots,p_n\in X$. Let
$L$ be a big and nef divisor on $X$ and let $F=(\sqrt{n/L^2})\pi^*L-E_1-\cdots-E_n$. 
Assume $H = \pi^* C -
(m+k)E_1-mE_2-\cdots-mE_n$ is the class 
of an almost uniform $F$-abnormal curve $H$ with
$n> m>0$. Then $-\sqrt{m} \le
k \le \sqrt{m}$. 
Moreover, if $k\ne0$, then also $C^2=2mk+m^2n$ (and
so $H^2=-k^2$) and $m\sqrt{n}-1 < \sqrt{C^2}\le C\cdot L/ \sqrt{L^2}
<m\sqrt{n}+1$.
\end{cor}

\eatit{
\begin{proof} We leave the case that $n=1$ to the reader, so assume $n\ge 2$.
To see $-\sqrt{m} \le k \le \sqrt{m}$,
observe that $m<n$ implies $mn/(n-1)\le m+1$; now apply
$k^2< mn/(n-1)$ from Corollary \ref{SzCor}(b).

Now, assume that $k\ne0$. By Corollary \ref{SzCor}(c) we have 
$(C\cdot L)^2/L^2-nm^2-2mk<k^2/n$, but now $k^2/n<1$; Corollary \ref{SzCor}(c) also
tells us that $C^2-nm^2-2mk\ge 0$. Therefore, 
putting both inequalities together with $C^2 \le (C\cdot L)^2/L^2$, we
must have $C^2-nm^2-2mk=0$, proving $H^2=-k^2$.

Finally, as $H$ is $F$-abnormal, we have
$(C\cdot L)/\sqrt{L^2}<m\sqrt{n}+k/\sqrt{n}<m\sqrt{n}+1$. 
On the other hand, since
$-k^2=H^2\ge -(m+k)$ by \cite{refX}, completing the square gives 
$(k-1/2)^2\le m + 1/4 < n$, so $k>1/2-\sqrt{n}$ and
$\sqrt{C^2}=\sqrt{m^2n+2mk}>\sqrt{m^2n-2m(\sqrt{n}-1/2)}\ge
\sqrt{(m\sqrt{n}-1)^2}$, 
and we conclude $\sqrt{C^2}>m\sqrt{n}-1$.
\end{proof}
}

\begin{Rmk}\label{lemrem} \rm
We note that if the N\'eron-Severi group of $X$ is generated by a
single ample  divisor $L$ with $L^2=r^2$ a square, when moreover
Corollary \ref{lemAAA} applies, there is for each $m$ at most one
$k\ne 0$ and one $t$ for which an abnormal curve
$[H] = t \pi^* L - (m+k)E_1-mE_2-\cdots-mE_n$ could exist.
Indeed, $t^2r^2=2mk+m^2n$ implies that $t^2r^2$ has the 
same parity as $m^2n$, and
only one integer $tr$ in the range 
$m\sqrt{n}-1<tr<m\sqrt{n}+1$ has this
property. 
\end{Rmk}

The next corollary is just
a refined version of Corollary \ref{SzCor}. Note that
$${\sqrt{L^2\over n}}\sqrt{1-{1\over{\mu n}}} 
= {\sqrt{L^2\over{n+\delta}}}$$
is equivalent to $\delta=(\mu-1/n)^{-1}$.
We will denote an almost uniform class of the form
$\pi^* C -m(E_1+\cdots+E_n)-kE_i$ by $H(C,m,k)$, with
$n$ being understood.

\begin{cor}\label{almunif}
Let $L$ be
a big and nef divisor on $X$. Let $\pi:Y\rightarrow X$ be the blow
up of $n>1$ general points $p_1,\ldots,p_n\in X$. Let $\mu\ge1$ be real
and consider the ${\bf R}$-divisor class 
$F(\delta)=\sqrt{(n+\delta)/L^2}L-(E_1+\cdots+E_n)$,
where $\delta=(\mu-1/n)^{-1}$. Then any $F(\delta)$-abnormal class is 
of the form $H(C,m,k)$,  where $C$, $m$ and $k$ 
are as in Corollary \ref{SzCor} and where $0<m<\mu$ and either $k=0$ or 
$m(n-1)<\mu$.
\end{cor}

\begin{proof} Let $H$ be an $F(\delta)$-abnormal class. Then $H=H(C,m,k)$, 
where $C$, $m$ and $k$ satisfy the criteria of Corollary \ref{SzCor}. 
First, say $k=0$; then $m^2n-m\le (C\cdot L)^2/L^2$, while
$F(\delta)\cdot C<0$ implies $(C\cdot L)\sqrt{(n+\delta)/L^2} < mn$, hence
$m^2n-m<m^2n^2/(n+\delta)$ 
or $(1/n)(1-1/(mn)) < 1/(n+\delta)$. This simplifies to
$m - 1/n<1/\delta=\mu - 1/n$, or $m<\mu$.
Now assume $k\ne 0$. This time we have $(C\cdot L)\sqrt{(n+\delta)/L^2} < mn+k$
and $m^2n+2mk+\hbox{max}\,(k^2-m,k^2-(m+k),0)\le (C\cdot L)^2/L^2$,
hence $(m^2n+2mk)/(mn+k)^2\le (C\cdot L)^2/((mn+k)^2 L^2)$. 
Note that $(1/n)(1-1/(mn(n-1))) \le (m^2n+2mk)/(mn+k)^2$
is the same as $1-1/(mn(n-1)) \le 
(m^2n^2+2mkn)/(mn+k)^2=1-k^2/(mn+k)^2$
or $mn(n-1)k^2\le (mn+k)^2$. This holds when $k>0$
because in this case $k^2<mn/(n-1)$. It also holds when $k<0$, because
now $k^2<(m+k)n/(n-1)$ or $mn(n-1)k^2<(m+k)mn^2$, but 
$(m+k)mn^2\le (mn+k)^2$ holds since it simplifies to
$kmn(n-2)< k^2$, but $k$ is negative. So, putting everything together, we have
$${1\over n}\left(1-{1\over{mn(n-1)}}\right) \le 
{{m^2n+2mk}\over{(mn+k)^2}} \le {{(C\cdot L)^2}\over{(mn+k)^2L^2}} < 
{1\over{n+\delta}}.$$ 
But $(1/n)(1-1/(mn(n-1))) < 1/(n+\delta)$ simplifies to
$m(n-1) - 1/n<1/\delta=\mu - 1/n$, or $m(n-1)<\mu$.\end{proof}

We can now prove Theorem \ref{thetheorem}, Corollary \ref{corAH} and Corollary \ref{simplifiedampCor}:
 
\begin{proof}[Proof of Theorem \ref{thetheorem}] Let us prove part (b)
  of Theorem \ref{thetheorem} first. 
Since $\sqrt{{L^2\over (n+\delta)}}=\sqrt{L^2\over n}\sqrt{1-{1\over \mu n}}$, 
the statement that \sc{X,L,n} is at least as big as 
$\sqrt{L^2\over n}\sqrt{1-{1\over \mu n}}$ follows if
$F(\delta)=\sqrt{(n+\delta)/L^2}L-(E_1+\cdots+E_n)$ is nef. 
If $F(\delta)$ were not nef, then there would exist an
$F(\delta)$-abnormal class $H=H(C,m,k)$, hence
$0> F(\delta)\cdot H$, so $(nm+k)/\sqrt{L^2/(n+\delta)} > L\cdot C \ge
\alpha_0((m^{[n-1]},m+k))$. But our hypotheses on
$\alpha_0$, together with Corollary \ref{SzCor} and Corollary \ref{almunif},
guarantee that this cannot happen.

Now consider (a). For every integer
$1 \le m < \mu$, assume that
$$\alpha(m^{[n]}) \ge m \sqrt{L^2(n-1/\mu)}
> m \sqrt{L^2(n-1/(\mu(1-2/(n+1))))}.$$
Then, whenever $1 \le m < \mu'=\mu(1-2/(n+1))$, 
we claim that $\alpha_0(m^{[n]}) \ge m \sqrt{L^2(n-1/\mu')}$,
and whenever $1 \le m < \mu'/(n-1)$,
$k^2 < (n/(n-1))\hbox{min}(m,m+k)$,
we claim that  $\alpha_0((m^{[n-1]},m+k)) \ge ((mn+k)/n)
\sqrt{L^2(n-1/\mu')}$. Part (b) will then imply that
$$\sc{X,L,n}\ge \sqrt{L^2/n}\sqrt{1-1/(n\mu')}>
\sqrt{L^2/n}\sqrt{1-1/((n-2)\mu)},$$ 
as wanted.

The first claim is immediate, for $m < \mu'< \mu$, so
$$\alpha_0(m^{[n]})\ge\alpha(m^{[n]}) \ge m \sqrt{L^2\left(n-{1\over
\mu}\right)}>m \sqrt{L^2\left(n-{1\over {\mu'}}\right)}.$$ 
For the second claim, given a reduced and 
irreducible curve $C=C_n$ with
multiplicity $m$ at general points $p_1$, \dots, $p_{n-1}$,
multiplicity $m+k$ at $p_n$ and $C\cdot  L=\alpha_0((m^{[n-1]},m+k))$, 
consider curves $C_1, \dots, C_{n-1}$
such that $C_i$ has multiplicity $m+k$ at $p_i$ and multiplicity $m$
at the other points (which exist because the points are general). 
Then $D=C_1+ \cdots + C_n$ is a (reducible) curve with
multiplicity $nm+k$ at each of the points.
But $k^2 < (n/(n-1))\hbox{min}\,(m,m+k)$ implies that
$k\le m$ (since otherwise $k^2\ge (m+1)^2 > 2m \ge nm/(n-1)\ge 
(n/(n-1))\hbox{min}(m,m+k)$, but this contradicts
Corollary \ref{SzCor}(b)). So if $m<\mu'/(n-1)$, then
$nm+k \le (n+1)m < (n+1)\mu'/(n-1)= \mu$, and
$$\hbox{$\alpha_0((m^{[n-1]},m+k)) \ge {1 \over n} \alpha((nm+k)^{[n]}) \ge
{{nm+k}\over{n}}\sqrt{L^2\left(n-{1\over{\mu'}}\right)}=
{{nm+k}\over{\sqrt{n}}}\sqrt{L^2\left(1-{1\over{n\mu'}}\right)},$}$$
as claimed.
\end{proof} 

\begin{proof}[Proof of Corollary \ref{corAH}]
Note that $\alpha(X,L,m^{[n]})/L^2$ is an integer which 
increases with $n$. Thus the Riemann-Roch formula together 
with ampleness of $L$ gives that 
$$\dim \Big|{{\alpha(X,L,m^{[n]})}\over{ L^2}}L\Big|=
{{\alpha(X,L,m^{[n]})(\alpha(X,L,m^{[n]})-L\cdot K)}\over{2 L^2}}+p_a
\,, $$
where $K$ denotes the canonical class and $p_a$ the arithmetic genus
of the surface $X$, provided that $\alpha(X,L,m^{[n]})$ is large
enough, which certainly holds (independent of
$m\ge 1$) for $n$ large enough.
Thus, in order to apply Theorem \ref{thetheorem}, it will be enough to prove 
for $n$ large enough that
$${{\alpha(\alpha-L\cdot K)}\over{2 L^2}} \ge n{{m(m+1)}\over 2} -p_a$$
implies $\alpha^2\ge m^2L^2 n$. If $L\cdot K\ge 0$ this is clear, so
assume $(L\cdot K)/L^2=-\beta<0$. Then, in order to have $\alpha^2<
m^2L^2 n$ it would be necessary that $\beta \alpha/2 > nm/2-p_a$ or
$\alpha>nm/\beta-c$ with $c=2p_a/\beta$ independent of $n$ and
$m$. But then $\alpha^2 >(nm/\beta
-c)^2 \ge m^2L^2n(n/(\beta^2L^2)-2c/(m\beta
L^2))\ge m^2L^2n(n/(\beta^2L^2)-2c/(\beta L^2))$, and 
for $n$ large enough this is 
bigger than $m^2L^2 n$, as desired. So it suffices to pick
$n_0$ large enough, then for $n\ge n_0$ we obtain the claimed lower bound
on $\sc{X,L,n}$.

We now verify that such an ${\bf m}(n)$ exists. 
Indeed, thanks to \cite{refAH}, a map ${\bf n}:{\bf N
\rightarrow N}$ exists such that for $n>{\bf n}(m)$ 
the inequality (\ref{ahbound}) holds.
Among such maps we may clearly choose one 
which is increasing. So, defining ${\bf m} : {\bf N\rightarrow N}$ 
as ${\bf m}(n)= \min\{m|{\bf n}(m)>n\}$, we have
for every $m< {\bf m}(n)$ that (\ref{ahbound}) holds.
Moreover, ${\bf m}$ is nondecreasing and unbounded 
since ${\bf n}$ is increasing,
hence $\displaystyle0=\lim_{n\to \infty}1/{\bf m}(n)=
\lim_{n\to \infty}n{\cal R}_n(L)$.
(Although \cite{refAH} does not give an explicit $\bf n$, we
have been informed by the authors that one may take ${\bf n}(m)\simeq
\exp(\exp(m))$, in which case ${\bf m}(n)\simeq \log(\log(n))$.)
\end{proof}

\begin{proof}[Proof of Corollary \ref{simplifiedampCor}]
Let $\delta= (\lfloor\sqrt{n}\rfloor(\lfloor \sqrt{n}\rfloor-3)/2 + 1-1/n)^{-1}$.
Then by Corollary \ref{corRem} and  the discussion immediately before  
Corollary \ref{almunif} we have $\sqrt{n+\delta}\ge \varepsilon(n)^{-1}$.
By hypothesis, $m<1/\sqrt{\delta}$, so $1/m^2>\delta$
so $\sqrt{n+1/m^2}>\sqrt{n+\delta}$.
Now $sL-E_1-\cdots-E_n$ is nef by Lemma \ref{introlem}(c),
for every rational $s$ such that $s\ge\sqrt{n+\delta}$,
hence $F(\delta)$ is itself nef.
Of course, $F\cdot E_i>0$ for all $i$. For any other reduced irreducible
curve $C$ it is enough to show
$C\cdot F(\delta)< C\cdot F$, since $0\le C\cdot F(\delta)$,
and $C\cdot F(\delta)< C\cdot F$ will follow if $t/m > \sqrt{n+\delta}$.
But $t^2\ge m^2n+1$, so $t/m\ge \sqrt{n+1/m^2}>\sqrt{n+\delta}$, as needed. 
\end{proof}

\end{document}